\documentclass[runningheads]{llncs}
\usepackage[T1]{fontenc}
\usepackage{amsmath}
\usepackage{amsfonts}
\usepackage{mathtools}
\usepackage{graphicx}
\usepackage{amssymb}
\usepackage{xcolor}
\usepackage{todonotes}

\mathtoolsset{showonlyrefs}
\newcommand{\R}{\mathbb{R}}

\newcommand{\p}{\mathcal{P}}

\newcommand{\XX}{\mathbb{X}}
\newcommand{\YY}{\mathbb{Y}}

\newcommand{\dx}{\,\mathrm{d}}
\newcommand{\eps}{\varepsilon}
\newcommand{\M}{\mathcal{M}}

\DeclareMathOperator{\OT}{OT}
\DeclareMathOperator{\UOT}{UOT}
\DeclareMathOperator{\GW}{GW}
\DeclareMathOperator{\W}{W}
\DeclareMathOperator{\FGW}{FGW}
\DeclareMathOperator{\UFGW}{UFGW}
\DeclareMathOperator{\UGW}{UGW}

\DeclareMathOperator{\KL}{KL}
\DeclareMathOperator{\id}{id}

\newcommand{\tT}{\mathrm{T}}

\begin{document}
\title{Gromov--Wasserstein Transfer Operators\thanks{Supported by the German Research Foundation (DFG) within the RTG
2433 DAEDALUS}
}
%
%
%
\authorrunning{
F. Beier
}
%
\author{Florian Beier
\orcidID{0000-0002-1692-9135},
\institute{
Institute of Mathematics, Technische Universit\"at Berlin, Stra\ss{}e des 17. Juni 136, 10623 Berlin, Germany
}
\email{f.beier@tu-berlin.de}
}
\maketitle              
\begin{abstract}
Gromov--Wasserstein (GW) transport is inherently invariant 
under isometric transformations of the data.
Having this property in mind, we propose to estimate dynamical systems 
by transfer operators derived from GW transport plans, 
when merely the initial and final states are known. 
We focus on entropy regularized GW transport, 
which allows to utilize the fast Sinkhorn algorithm and 
a spectral clustering procedure to extract coherent structures. 
Moreover, the GW framework provides a natural quantitative assessment 
on the shape-coherence of the extracted structures.
We discuss fused and unbalanced variants of GW transport for labelled and noisy data, respectively. 
Our models are verified by three numerical examples of dynamical systems with governing rotational forces.

\end{abstract}
%
%
%

\section{Introduction}
Optimal transport (OT) aims to find an optimal mass transport between two input (marginal) measures 
according to an underlying cost function. To improve the speed of the numerical computation,
Cuturi \cite{cuturi} introduced a regularized OT version which can be solved by the
fast and parallelizable Sinkhorn algorithm. 
Further effort has been made to generalize the OT for different settings 
as, e.g., unbalanced optimal transport \cite{LMS18}, 
which relaxes the hard matching of the marginal measures.
Another line of work pioneered by Mémoli \cite{memoli} focuses 
on so-called Gromov--Wasserstein (GW) distances. 
Here, the inputs have additional structure in the sense of intrinsic (dis-)similarities. 
The difference to $\OT$ is that a meaningful cost function on the product space of the inputs might not be available. 
Instead, the mass is transported so that pairwise (dis-)similarities are preserved. 
GW distances are invariant under isometric transformations, 
making them a valuable tool for e.g.\ shape classification \cite{lgw}, word alignment \cite{gw_words} or graph matching \cite{gw_graphs}.
For certain applications, a transport which simultaneously takes 
structural data in the GW sense as well as labelled data in the OT sense into account,
is desirable. This is possible in the framework of fused GW transport \cite{fgw}. 
Moreover, (fused) GW transport allows for a similar entropic regularization and unbalanced relaxation as OT \cite{umgw,ugw,ufgw}. 
In \cite{invOT} the authors propose a framework which extends OT to be invariant to various classes of linear transformations such as e.g.\ orthogonal transformations. Compared to GW this method is numerically more appealing but has the drawback that it requires the inputs to be embedded in a common space and centered.

Recently, Koltai et al. \cite{KLNS20} examined OT-based estimations of dynamical systems 
from observed initial and final states.
More precisely, the authors 
leveraged solutions to regularized (unbalanced) optimal transport 
to estimate so-called transfer operators. 
These are linear operators that characterize dynamical systems in the form of density flows. 
Furthermore, a clustering procedure based on the spectral information of the estimator 
was used to extract so-called coherent structures of the dynamical system. 
Although no unified definition of such structures is available, 
it is understood that they are persistent in time and space. 
Coherent structures are of particular interest e.g.\ in fluid dynamics, since they capture important flow dynamics.
This makes precise knowledge of the formation of coherent structures very appealing,
since it may lead to a deeper understanding of the dynamics or computational advancements.
In \cite{junge2022entropic}, the authors assumed instead
that the exact transfer operator is known on a finite subset of the full state space. 
Then, using regularized OT, a finite-dimensional approximation is constructed which limit 
is a regularized version of the ground truth and exhibits desirable properties, 
such as retention of the spectral information. 

In this paper, we build on the work  in \cite{KLNS20}, 
but use entropic GW transport plans for constructing transfer operators. 
This is motivated by the fact that GW transport 
is readily able to detect isometric transformations such as rotation. 
Additionally, data labels can be incorporated.
We will see that our proposed model includes a quantitative assessment of shape-coherence of the extracted structures.

{\bf Outline of the paper.}
In Section \ref{sec:prel}, we briefly recall regularized (unbalanced) OT,
associated transfer operators and related spectral clustering procedures.
GW transport and its (unbalanced) regularized and fused variants are introduced 
in Section \ref{sec:gw_trans}.
Then, we expand the derivation of transfer operators and spectral clustering towards GW transport plans.
In Section \ref{sec:num}, we present numerical examples which indicate the potential of our method.
%
\section{Optimal Transport and Transfer Operators}\label{sec:prel}
We consider (unbalanced) entropic OT, show how
transfer operators can be derived from OT plans, and elaborate on spectral clustering.
The derivation of transfer operators will be generalized to GW plans
in the next section.
\\[1ex]
{\bf Optimal transport.}
Let $X,Y \subset \R^d$ be compact sets equipped with the Euclidean distance $d_\text{E}$. 
By $\M^+(X)$, we denote the set of non-negative (Borel) measures and by $\p(X) \subset \M^+(X)$ the set of probability measures on $X$. 
Furthermore, let $L_{\mu}^2(X)$ be the Hilbert space of (equivalence classes) of square integrable functions with respect to the finite measure 
$\mu \in \M^+(X)$
equipped with the inner product $\langle \cdot,\cdot \rangle_\mu$. 
By $1_A$, we denote the characteristic function on $A$.
For $\mu,\nu\in {\mathcal M^+}(X)$, the \emph{Kullback--Leibler divergence}  is defined by
\begin{equation*}
\KL (\mu,\nu) \coloneqq \int_{X} \log\bigl(\tfrac{\dx \mu}{\dx \nu}\bigr) \dx \mu + \nu(X) - \mu(X),
\end{equation*}
if the Radon--Nikodym derivative $\frac{\dx \mu}{\dx \nu}$ exists,
and by $\KL(\mu,\nu) \coloneqq \infty$ otherwise. 
For $\mu \in \p(X)$, $\nu \in \p(Y)$, a lower semi-continuous cost function $c: X \times Y \to [0,\infty)$ and $\eps > 0$,
the \emph{regularized OT problem} is given by
\begin{equation} \label{ot_e}
\OT_\eps(\mu,\nu) \coloneqq \min_{\pi \in \Pi(\mu,\nu)} \underbrace{\int_{X \times Y} c(x,y) \dx \pi + \eps \KL(\pi, \mu \otimes \nu)}_{\eqqcolon F_\varepsilon^{\OT}(\pi)},
\end{equation}
where $\Pi(\mu,\nu) \coloneqq \{\pi \in \p(X \times Y) : {P_1}_\# \pi = \mu, {P_2}_\# \pi = \nu\}$ with
$P_i(x_1,x_2) \coloneqq x_i$ and push forward measures ${P_i}_\# \pi = \pi \circ P_i^{-1}$, $i=1,2$. Elements of $\Pi(\mu,\nu)$ are called \emph{transport plans}.
For $\eps = 0$, we obtain the unregularized optimal transport $\OT(\mu,\nu)$. 
The minimizer in \eqref{ot_e} is called (entropic) optimal transport plan $\hat \pi_\varepsilon$.
In the following, we will mainly use $c = d_{\text{E}}^2$ which leads to the \emph{Wasserstein distance} $\OT(\mu,\nu)^\frac12$.
The dual problem of $\OT_\varepsilon$ is
	\begin{align*}
	\OT_\varepsilon(\mu,\nu) 
	&= \max_{{\tiny (f,g)\in L^\infty_\mu(X)\times L^\infty_\nu(Y)}}
	\Big\{ \int_{X}f \dx \mu + \int_{Y} g \dx \nu \\
	&\quad - \varepsilon \int_{X \times Y} \exp\Bigl(\frac{ f(x) + g(y) - c(x,y)}{\varepsilon}\Bigr) -1
	\dx (\mu \otimes \nu) \Big\}
		\end{align*}
	Optimal potentials $\hat f_\varepsilon \in L^\infty_\mu(X)$, $\hat g_\varepsilon \in L^\infty_\nu(Y)$ 
	exist and are unique on ${\rm supp}(\mu)$ and ${\rm supp}(\nu)$ up to an additive constant.
	They are related to $\hat \pi_\varepsilon$ by
	\begin{equation} \label{eq:PDrelation}
	\hat\pi_\varepsilon = \exp\Bigl(\frac{\hat f_\varepsilon(x) + \hat g_\varepsilon(y) - c(x,y)}{\varepsilon}\Bigr) \, (\mu \otimes \nu)
	\eqqcolon  k_{\varepsilon} (\mu \otimes \nu).
	\end{equation}
	For atomic measures, the solution can be approximated efficiently by Sinkhorn's algorithm.
In some applications, it is useful to deal with \emph{regularized unbalanced OT}
\begin{align*}
\UOT_{\eps,\kappa} (\mu,\nu) 
\coloneqq  \min_{\pi \in \M^+(X\times Y)} F_\eps^{\OT}(\pi)
+ \kappa \left( \KL({P_1}_\# \pi, \mu) + \KL({P_2}_\# \pi, \nu) \right), \, \kappa > 0,
\end{align*}
which relaxes the hard marginal constraints on the objective to penalizing the $\KL$ divergence of its marginals with respect to the inputs. Unbalanced optimal transport is treated in detail in \cite{LMS18} and its regularized version in \cite{SFVTP19}.
Similarly as in the balanced case, there is a dual problem formulation 
with optimal potentials $(\hat f_{\varepsilon,\kappa}, \hat g_{\varepsilon,\kappa}) \in L^\infty_\mu(X) \times L^\infty_\nu(Y)$ 
and the optimal transport plan is given by 
\begin{equation} \label{eq:ot-unbalanced}
\hat \pi_{\varepsilon,\kappa} = \exp\Bigl(\frac{\hat f_{\varepsilon,\kappa}(x) + \hat g_{\varepsilon,\kappa}(y) -c(x,y)}{\varepsilon}\Bigr) \, (\mu \otimes \nu)
\eqqcolon  k_{\varepsilon,\kappa} \,  (\mu \otimes \nu).
\end{equation}
A generalization of Sinkhorn's algorithm can be used to solve the corresponding discrete problem, see \cite{SFVTP19}.
\\[1ex]
{\bf Transfer operators.}
\emph{Transfer operators}, also known as \emph{Perron-Frobenius operators}, 
are linear operators which characterize dynamical systems in the form of density flows \cite{froyland}. 
We consider transfer operators derived from entropic transport plans as in \cite{KLNS20}. 
Here we restrict ourselves to the balanced setting as the unbalanced case follows in a similar way.
To this end,  we associate $\hat \pi_\eps$ in \eqref{eq:PDrelation} with the transfer operator
$K_\eps:L^2_\mu(X) \to L^2_\nu(Y)$ given by
\[
(K_\eps \psi)(y) \coloneqq \int_{X} k_\eps(x,y) \psi(x) \dx \mu(x). 
\]
Figuratively, $K_\eps$ captures the structure of the transport of $\hat{\pi}_\eps$ independent of the marginal masses. 
Since $\hat \pi_\varepsilon \in \Pi(\mu,\nu)$, it holds
\begin{align*}
\int_X k_\varepsilon (x,y) \dx \mu(x) = 1_Y \; 
\,\nu\text{-a.e.}
\quad \text{and} \quad
\int_Y k_\varepsilon (x,y) \dx \nu(y) = 1_X \; \,\mu\text{-a.e.}.
\end{align*}
In particular, for atomic measures 
$\mu = \sum_{i=1}^m \mu(i) \delta_{x_i}$ and $\nu = \sum_{j=1}^n \nu(j) \delta_{y_j}$ 
and an optimal transport plan 
$\hat \pi_\varepsilon = \sum_{i,j=1}^{m,n} \hat \pi_\varepsilon (i,j) \delta_{x_i,y_j}$ 
using the matrix-vector notation 
${\bf \mu} \coloneqq (\mu(i))_{i=1}^m$, 
$\nu \coloneqq (\nu(j))_{j=1}^n$, $D_\mu \coloneqq {\rm diag} (\mu)$ and
$\hat \pi_\varepsilon \coloneqq (\hat \pi_\varepsilon (i,j))_{i,j=1}^{m,n}$,
the transfer kernel and operator are given by
\[
k_\varepsilon =  D_\mu^{-1} \, \hat \pi_\varepsilon \, D_\nu^{-1}
\quad \text{and} \quad
K_\varepsilon = D_\nu^{-1} \, \hat \pi_\varepsilon^\tT.
\]
{\bf Spectral clustering.}
In order to extract coherent structures in dynamical systems, we can apply a spectral clustering procedure on $K_\varepsilon$, see \cite{KLNS20}.
The clustering premise is just the knowledge of two observations from $\mu \in \p(X)$ and $\nu \in \p(Y)$ in a dynamical system
without any knowledge of the true dynamics.
The goal is to find measurable partitions $X = X_1 \dot \cup X_2$, $Y = Y_1 \dot \cup Y_2$ fulfilling ideally
\begin{equation}\label{eq:coherence}
K_\eps 1_{X_k} = 1_{Y_k} \quad \text{and} \quad \mu(X_k) = \nu(Y_k), \quad k =1,2.
\end{equation}
These conditions may be interpreted as coherence and mass preservation of the partitions.
One way to tackle this problem is to consider the following optimization problem
\[
\max_{X_1 \dot \cup X_2 
= 
X,Y_1 \dot\cup Y_2 = Y}  \biggl\{ \frac{\langle K_\eps 1_{X_1}, 1_{Y_1} \rangle_\nu}{\mu(X_1)} + \frac{\langle K_\eps 1_{X_2}, 1_{Y_2} \rangle_\nu}{\mu(X_2)} \biggr\},
\]
which is usually relaxed to
\begin{equation}\label{eq:fuzzy}
\max_{(\varphi,\psi) \in L^2_\mu(X) \times L^2_\nu(Y)} \biggl\{\frac{\langle K_\eps \varphi,\psi\rangle_\nu}{\|\varphi\|_{\mu} \|\psi\|_{\nu}}: 
\langle \varphi, 1_X \rangle_\mu = \langle \psi,1_Y\rangle_\nu = 0\biggr\}.
\end{equation}
Since $K_\eps$ is bounded and non-negative $(\mu \otimes \nu)$-a.s., it follows that the largest singular value of $K_\eps^* K_\eps$ is simple \cite[Lem.~3]{froyland}. Moreover, the largest singular value of $K_\eps$ is $1$ and the corresponding left and right singular functions are $1_X$ and $1_Y$, respectively. Notably, $(1_X,1_Y)$ are not included by the constraints in \eqref{eq:fuzzy}.
Hence, a maximizing pair $(\hat{\varphi},\hat{\psi})$ in \eqref{eq:fuzzy} is given by the right and left singular functions of $K_\eps$ 
associated to the second largest singular value of $K_\eps$. The desired partitioning is then readily obtained by thresholding $(\hat{\varphi},\hat{\psi})$ at zero. Solving \eqref{eq:fuzzy} in practice amounts to computing a (truncated) singular
value decomposition of $K_\eps$.

\section{Transfer Operators from GW Transport Plans}\label{sec:gw_trans}
In both references \cite{junge2022entropic,KLNS20}, the assumption on the underlying dynamics is that they are compliant with an optimal transport. 
For certain situations this might not be the case. 
Consider e.g.\ particles on the two-dimensional unit disk with a driving rotational force. 
If the rotation angle between two observations is large, OT will not be able to recover this dynamic, see the first Example in Section \ref{sec:num}.
A transport setting which naturally handles isometric transforms such as rotation is given by the framework of GW transport \cite{memoli}. 
As before, we consider compact state spaces $X,Y \subset \R^d$ and measures $\mu \in \p(X)$, $\nu \in \p(Y)$.
In the contrast to classic OT, a cost function on the product space $X \times Y$ is not required. 
Instead we seek the preservation of the internal structure of the spaces. Here we focus on the Euclidean metrics,
for generalizations see \cite{sturm}.
We set $d_X \coloneqq d_\text{E} \vert_{X \times X}$.
Then the triples $\XX \coloneqq (X,d_X,\mu)$, $\YY \coloneqq (Y,d_Y,\nu)$ are called \emph{metric measure (mm-) spaces}. 
We introduce the notation $\mu^{\otimes} \coloneqq \mu \otimes \mu$.
For $\eps > 0$, the \emph{regularized GW transport} between two mm-spaces $\XX$ and $\YY$ is defined by
\begin{align}\label{eq:GW}
    &\GW_\eps(\XX,\YY) \coloneqq \inf_{\pi \in \Pi(\mu,\nu)} F_\varepsilon^{\GW}(\pi), \\
    &F_\varepsilon^{\GW}(\pi) \coloneqq \hspace{-0.3cm}\int\limits_{(X \times Y)^2} \hspace{-0.25cm}(d_X(x,x') - d_Y(y,y'))^2 \dx \pi(x,y) \dx \pi(x',y') + \eps \KL\bigl(\pi^{\otimes}, (\mu \otimes \nu)^{\otimes}\bigr).
\end{align}
In contrast to $\OT$, we regularize with the quadratic $\KL$ divergence as in \cite{ugw}. For $\eps = 0$, we obtain the unregularized GW transport $\GW$ which was originally introduced in \cite{memoli}.
Notably, $\GW^\frac12$ defines a metric on the space of mm-spaces up to identification by measure-preserving isometries. More precisely, $\GW(\XX,\YY) = 0$ if and only if there exists an isometry $I:X \to Y$ with $\nu = I_\# \mu$. In this case, $(\id,I)_\# \mu$ is an optimal GW plan. In particular, this shows the invariance of GW with respect to isometric transformations. Figuratively, optimal $\GW$ plans are such that whenever they transport (infinitesemal) mass from $x$ to $y$ and $x'$ to $y'$ one has $d_X(x,x') \approx d_Y(y,y')$ which favors a near-isometric transport.

Similar to $\OT_\eps$, $\GW_\eps$ admits unbalanced versions \cite{ugw}, we focus on marginal penalization using $\KL$. For $\eps,\kappa> 0$, the \emph{unbalanced regularized GW transport} is defined by
\begin{align*}\label{eq:UGW}
    \UGW_{\eps,\kappa}(\XX,\YY) = & \hspace{-0.35cm}\inf_{\pi \in \M^+(X \times Y)} \hspace{-0.15cm} F_\varepsilon^{\GW}\!(\pi) + \kappa \left( \KL(({P_1}_\# \pi)^{\otimes},\!\mu^{\otimes}) + \KL(({P_2}_\# \pi)^{\otimes}\!, \nu^{\otimes})  \right)\!.
\end{align*}
Here the marginals of optimal plans differ from the inputs whenever an exact matching results in large values under the functional $F_\eps^{\GW}$. This can make $\UGW_{\eps,\kappa}$ somewhat robust to outliers.

When working with labelled data, we might be interested in a transport plan which preserves the internal geometrical information in the form of metrics as well as feature information in the form of labels. This leads to a fused version of the GW and the Wasserstein distance. To incorporate label information, we introduce an additional set $A \subset \R^m$ endowed with $d_A \coloneqq d_{\text{E}}\vert_{A \times A}$. We assume that each point in $X,Y$ admits only one label, 
which we characterize by label functions $l_X:X \to A$, $l_Y:Y \to A$, respectively. 
Clearly, a more general treatment would be to consider distributions in the label space as in e.g.\ \cite{fgw}. 
In our case, the \emph{regularized fused GW distance} is defined by
\begin{equation*}
    \FGW_\eps((\XX,l_X),(\YY,l_Y)) \coloneqq \inf_{\pi \in \Pi(\mu,\nu)} F_\eps^{\GW}(\pi) + \int_{X \times Y} d_A\bigl(l_X(x),l_Y(y)\bigr) \dx \pi(x,y).
\end{equation*}
As with the original formulation, the marginal constraints may be relaxed in the same way which leads to an unbalanced, fused variant $\UFGW_\eps^\kappa$ which was discussed in \cite{umgw,ufgw}.

The previously discussed GW formulations are quadratic with respect to the objective plan which renders them numerically challenging.
For our numerical experiments below, we rely on a class of simple iterative algorithms which are based on block-coordinate relaxations. The main idea consists of alternately fixing one plan while minimizing with respect to the other. The problem that is then minimized in each iteration step can be written as an entropic OT problem for which Sinkhorn's algorithm can be leveraged.
Details regarding this procedure can be found in
\cite{peyreGW} (balanced GW), \cite{ugw} (unbalanced GW) and \cite{umgw,ufgw} (unbalanced, fused GW).
Solutions $\hat{\pi}_\eps$ obtained with this procedure are also solutions to an entropic (unbalanced) OT problem 
and thus have the form \eqref{eq:PDrelation}, i.e.\ it holds $\hat{\pi}_\eps = k_\eps \, (\mu \otimes \nu)$. 
Ultimately, this allows us to apply the spectral clustering procedure on the associated transfer operator $K_\eps$ as described in Section \ref{sec:prel}.
The next remark highlights another benefit of GW over OT transfer operators for extracting coherent structures.

\begin{remark}[Quantitative assessment of shape-coherence]\label{rem:1}
      Let $\hat \pi_\eps$ be an optimal GW plan between $\XX$ and $\YY$ with associated transfer operator $K_\eps$ and
			$X_i,Y_i$, $i=1,2$ the spectral clustering partition.
			Even if the partitions satisfy $K_\eps 1_{X_i} \approx 1_{Y_i}$ and $\mu(X_i) \approx \nu(Y_i)$, 
			it may be that the intrinsic shapes of $X_i$ and $Y_i$, $i = 1,2$ differs significantly. 
			It depends on the application, if these structures should be considered coherent or not. 
			The GW framework readily gives us the possibility for a quantitative assessment of shape-coherence by 
			evaluating the $\GW$ functional at $\hat \pi_\eps$ restricted to $X_i \times Y_i$, $i =1,2$. 
			The closer the evaluation is to $0$, the more the associated partitions can be considered shape-coherent 
			or isometric under the transfer operator $K_\eps$. We apply this for Example 3 in Section \ref{sec:num}.
\end{remark}

\section{Numerical Examples}\label{sec:num}

In this section we provide three examples of our proposed GW transfer method. 
In OT comparisons we use the quadratic Euclidean cost function.
We partly rely on the Python Optimal Transport library \cite{pot}. For our experiments we aim to set the entropic regularization parameter $\eps > 0$ as small as possible while avoiding numerical overflow.
\\[1ex]
{\bf 1. Particles on a rotating disk.}
First, we are interested in the ability of OT plans to recover the dynamics of a rotating system and compare it with a GW based approach.
\begin{figure}[t]
    \centering
    \includegraphics[width = 0.66\linewidth]{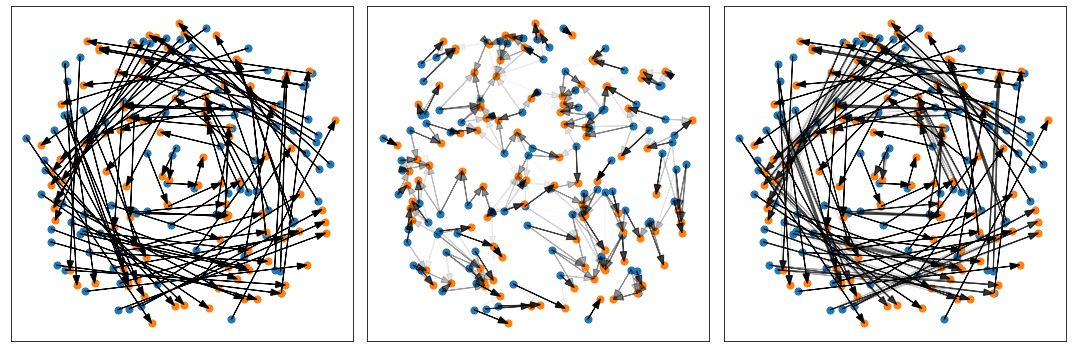}
    \caption{An example of $X$ (blue), $R_{\frac{\pi}{2}}(X)$ (orange), as well as the associated transfer kernels $k_\text{true}$ (left), $k^{\W}_\eps$ (middle) and $k^{\GW}_\eps$ (right). The arrows are drawn from $x$ to $y$ with opacity proportional to the respective kernel at $(x,y)$.}
    \label{fig:num1}
\end{figure}
We consider $n=50$ uniformly sampled particles on the 2D unit disk $D \coloneqq \{x : \|x\|_2 \leq 1\} \subset \R^2$. 
Let $X \subset D$ be the set of particles. 
We consider a counter-clockwise rotation of degree $\theta \in (0,2\pi)$. 
More precisely, the true transfer operator is characterized by the bijective map $R_\theta:D \to D$ given by
\[
(r \cos(\phi), r \sin(\phi)) \mapsto (r \cos(\phi + \theta), r \sin(\phi + \theta)), \quad r \in [0,1], \theta \in [0,2\pi).
\]
We focus on the transfer associated to  kernel $k_{\text{true}}(x,y) = \delta_{\{R_\theta(x) = y\}}$.
An illustration of $X$ and $R_\theta$ is shown in Figure \ref{fig:num1}.
We investigate how well the GW transfer operator estimates the true transfer operator 
for $\theta = \frac{\pi}{30}, \frac{2\pi}{30}, \dotsc, \pi$.
To this end, we sample the initial state $X$ and  compute the GW transport plan with $\eps = 0.0008$
between between the uniform distributions on $X$ and $R_\theta(X)$, respectively. 
As discussed, all plans admit the form \eqref{eq:PDrelation} for respective kernels $k_\eps^{\W}, k_\eps^{\GW}$.
To compare the performance we consider the error measure 
$$\text{e}(k_\varepsilon^\bullet) \coloneqq \frac{1}{n^2} \sum_{x \in X} \sum_{y \in R_\theta(X)} k_\varepsilon^\bullet  (x,y) d_{\text{E}}(R_{\theta}(x),y),
\quad \bullet \in \{{\rm W}, {\rm GW}\}.
$$
Intuitively, this gives us the mean Euclidean distance when comparing the transfer operator associated to the kernel against the true transfer.

\begin{figure}[t]
    \centering
    \includegraphics[width = 0.65\linewidth]{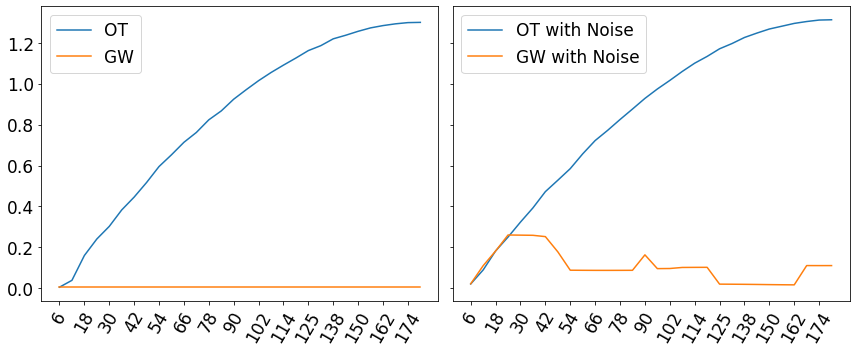}
    \caption{The mean errors plotted against rotation angle $\theta$ in degrees without noise (left) and with noise (right).}
    \label{fig:num1-2}
\end{figure}

\begin{figure}[t]
    \centering
    \includegraphics[width=\linewidth]{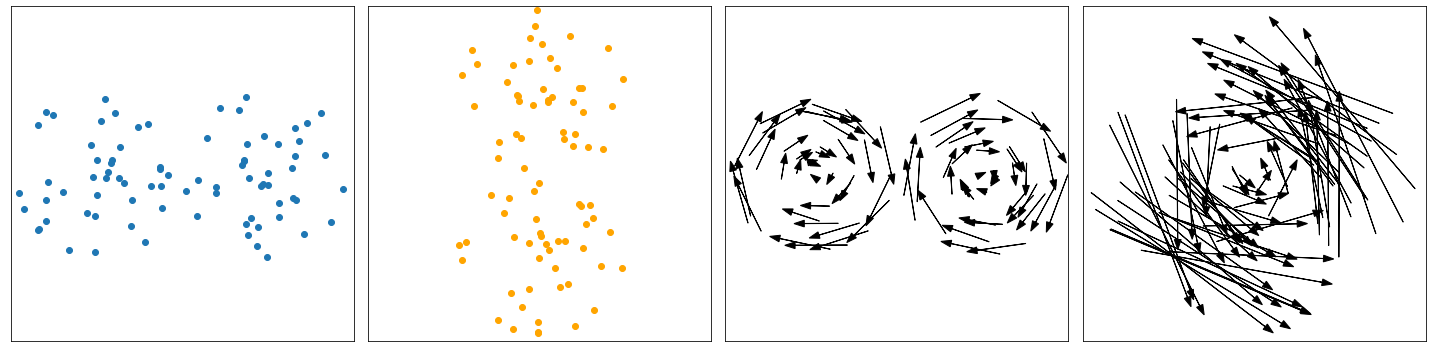}
    \caption{Left to right: $X$, $Y$, $(R^{(1)} + R^{(2)})$ and $R_{\theta}$.}
    \label{fig:num2_1}
\end{figure}

\begin{figure}[t]
    \centering
    \includegraphics[width = \linewidth]{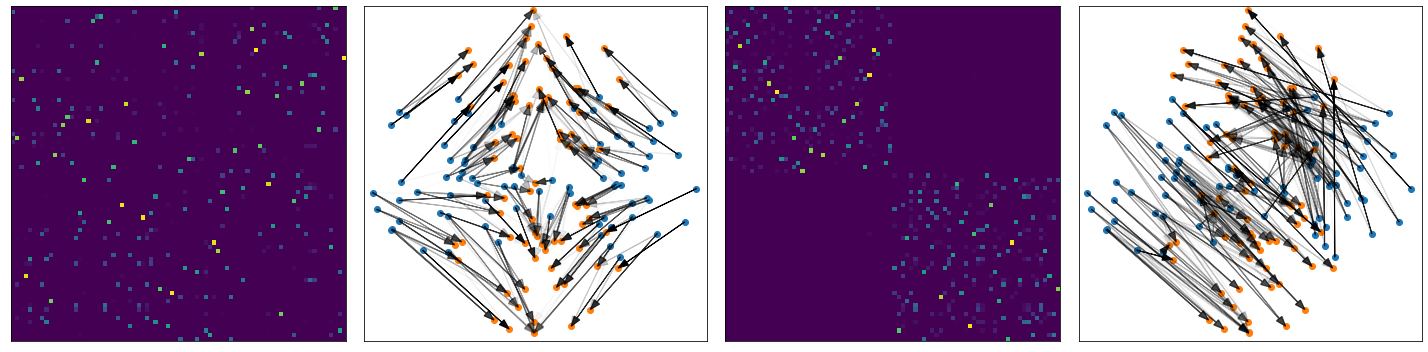}
    \caption{Left to right: Matrix $\pi^{\W}_\eps$, transfer kernel $k_\eps^{\W}$, matrix $\pi^{\GW}_\eps$, transfer kernel $k_\eps^{\GW}$.}
    \label{fig:num2_2}
\end{figure}

\begin{figure}[t]
    \centering
    \includegraphics[width = \linewidth]{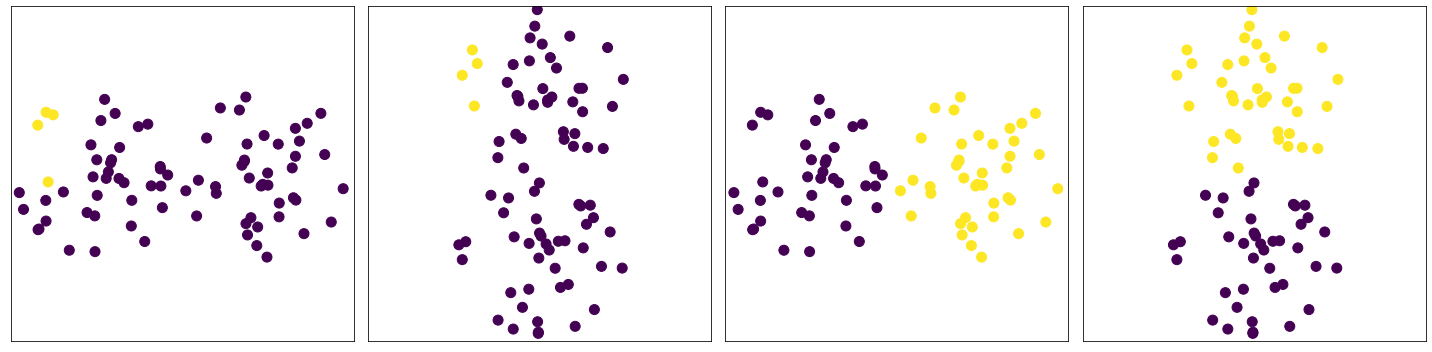}
    \caption{The spaces $X,Y$ coloured according to the sign of the left and right eigenvector of $K^{\W}_\eps$ (left) and $K^{\GW}_\eps$ (right) corresponding to the second largest eigenvalue.}
    \label{fig:num2_3}
\end{figure}

The right-hand side of Figure \ref{fig:num1} shows the qualitative difference between the OT and GW-based approaches for one example with $\theta = \frac{\pi}{2}$.
A quantitative comparison is given on the left-hand side of Figure \ref{fig:num1-2}. 
More precisely, we sampled $10$ independent choices of $X$ to 
obtain $10$ OT plans $\pi_\eps^{\W,1},\dotsc,\pi_\eps^{\W,10}$ and $10$ GW plans $\pi_\eps^{\GW,1},\dotsc,\pi_{\eps}^{\GW,10}$ for each angle $\theta$. 
We plot the mean errors $\frac{1}{10}\sum_{i=1}^{10} \text{e}(k_\eps^{\bullet,i})$, $\bullet \in \{{\rm W}, {\rm GW}\}$
as a function of the angle $\theta$.
As expected, for large values of $\theta$, the OT-based transfer operator is a poor estimator. 
This is evident since, e.g.\ for a 90 degrees rotation, points are transferred far distances which is sub-optimal in the OT sense. 
Even for smaller angles such as 18 degrees, we observe a mean error of $0.15$. 
On the other hand, the GW based approach recovers $R_\theta$ nearly exactly in all cases.
\begin{figure}[t]
    \centering
    \includegraphics[width = 0.6\linewidth]{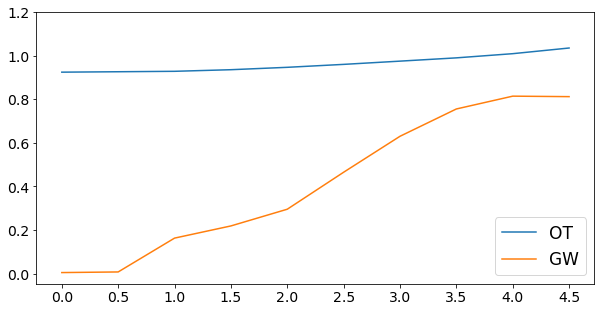}
    \caption{The mean errors for fixed $\theta = \frac{\pi}{2}$ plotted against the noise magnitude $m$.}
    \label{fig:num1-3}
\end{figure}
In the previous example, $Y = R_\theta(X)$ was given by an exact rotation of $X$. 
However, in practice the observed end state $Y$ of the dynamical system might be a noisy version of $R_\theta(X)$.
Hence, we repeat the previous experiment, where this time $Y = R_\theta(X) + m \eta$ with $\eta \sim \mathcal U([-0.1,0.1]^2)$ and $m=1$. 
To make this comparable to the previous experiment, we consider the same sampled initial states $X$ as above. 
We proceed as before and plot the error of the OT and GW-based approach on the right-hand side of Figure \ref{fig:num1-2}. 
For small angles, $\GW$ remains comparable to $\OT$ whereas for large angles a better estimation is achieved by utilizing $\GW$.
Finally, we repeat the procedure this time for a fixed angle $\theta = \frac{\pi}{2}$ and for varying $m = 0.5,1,1.5,\dotsc,4.5$. The result is plotted in Figure \ref{fig:num1-3}.
\\[1ex]
{\bf 2. Multiple rotating disks.}
In our next example, let $\theta = \pi/2$, and $D$, $R_{\theta}$ as above. In addition, for $i =1,2$, we consider
\[
D_i = \{x \in D : \|x - x^{(i)}\| \leq 1/2\}, \quad x^{(i)} = (-1/2,0), \quad x^{(i)} = (1/2,0).
\]
We set $F = (R^{(1)} + R^{(2)}) \circ R_{\theta}$, where $R^{(i)}$ constitutes a rotation of $- \pi / 4$ around $x^{(i)}$, restricted to $D_i$, $i=1,2$. 
Let $n = 80$, we uniformly sample $n/2$ points of $D_1$ and $D_2$, respectively. Denote the entire set of $n$ points by $X \subset D$. 
Let $Y = F(X)$ and equip $X$ and $Y$ with the uniform distribution. Figure \ref{fig:num2_1} illustrates $X$, $Y$ and $F$.
We focus on the estimation of the transfer operator associated to $k_{\text{true}}(x,y) = \delta_{\{F(x) = y\}}$.
We  compute an OT plan $\pi^{\W}_\eps$ and an GW plan $\pi^{\GW}_\eps$ both with $\eps = 0.001$. 
Illustrations of the matrices $\pi^{\W}_\eps,\pi^{\GW}_\eps$ as well as a visualizations of the transfer operators $K^{\bullet}_\eps$, associated to respective kernels $k^{\bullet}_\eps$, $\bullet \in \{{\rm W}, {\rm GW}\}$
are provided in Figure \ref{fig:num2_2}. 
Clearly, neither approach is able to recover the ground truth.
However, the figure indicates that $K^{\GW}_\eps$ transfers most of the mass from $D_i$ to $F(D_i)$, $i =1,2$, while the OT-based approach does not. 
We apply the spectral clustering procedure, i.e. we compute the left and right eigenvectors associated to the second largest eigenvalue of 
$K^{\bullet}_\eps$, $\bullet \in \{{\rm W}, {\rm GW}\}$
and present them in Figure \ref{fig:num2_3}.
As expected, the partitioning according to $K^{\GW}_\eps$ is able to find both coherent disks. 
Now, the local dynamics within the partitions can readily be obtained by computing the GW transport of the partitioned subspaces.

\begin{figure}[t]
    \centering
    \includegraphics[width = \linewidth]{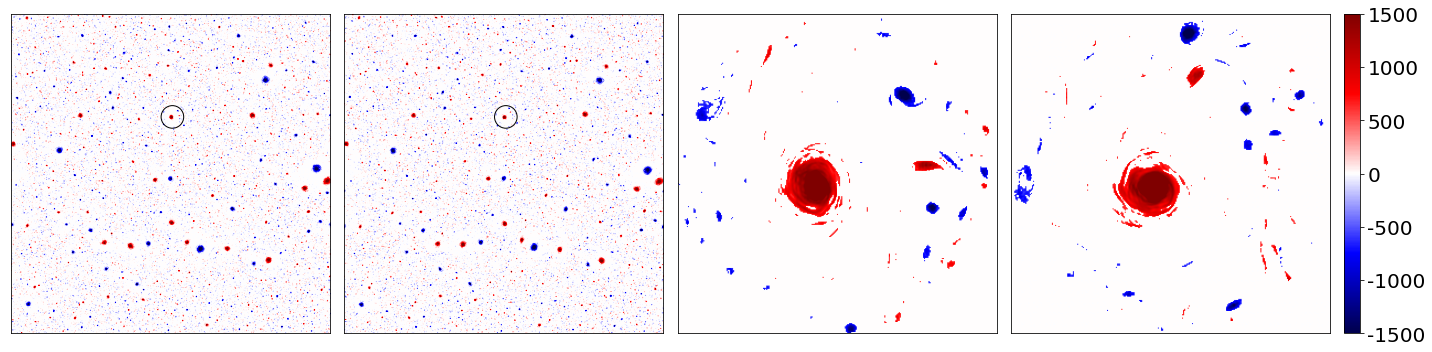}
    \caption{Left: Two thresholded snapshots $\omega_0,\omega_1$ of a direct numerical simulation of \eqref{eq:w} on a 4096x4096 pixel grid. The images on the right-hand side show the marked circular patch for both time-steps.}
    \label{fig:five_snaps}
\end{figure}

We conclude this example by remarking that the correct identification of the discs may also fail and is not stable with respect to noise.
This is due to the fact that e.g.\ $\XX$ is almost isometric to a 180 degrees rotation as well as a reflection along the vertical axis. If the inputs are subjected to noise, an optimal $\GW$ plan might match $D_1$ with $F(D_2)$ and $D_2$ with $F(D_1)$.
\\[1ex]
{\bf 3.  Vorticity field of the 2D Navier--Stokes equation.}
Finally, we consider a two-dimensional flow in time which behaves according to the 2D Navier--Stokes equations on the square $[0,2\pi]^2$ (periodic boundary conditions)
\begin{align*}
    \partial_t u + (u \cdot \nabla) u &= - \nabla p + v \nabla^2u\\
    \nabla \cdot u &= 0,
\end{align*}
where 
$u: [0,T] \times [0,2\pi]^2 \to \R^2$ is the velocity, $p:[0,T] \times [0,2\pi]^2 \to \R^2$ the pressure and $v \in \R$ the kinematic viscosity. 
Numerically, it is more efficient to solve the 
scalar advection-diffusion equation
\begin{equation}\label{eq:w}
\partial_t \omega + (u \cdot \nabla) \omega = v \nabla^2 \omega,
\end{equation}
where $\omega = \partial_x u_y - \partial_y u_x$ is the vorticity of $u$. 
Following \cite[Sec~IV]{wang2022coherent}, the equation is solved in the Fourier domain after a adding a small-scale forcing term 
and a large-scale damping function on a 4096x4096 grid.
Ultimately, we obtain two time snapshots $\omega_0,\omega_1$ of the vorticity field, which we restrict to $\lvert \omega_i \rvert \geq 600$. 
The snapshots as well as a zoom into a circular patch with a $290$ pixel diameter is shown in Figure \ref{fig:five_snaps}.
As we can see, the flow exhibits coherent structures in the form of vortices on large and small scales. 
Large vortices essentially determine most of the local dynamics. This can be seen for instance in the in selected patch, where smaller vortices are rotating around the large center vortex.
\begin{figure}[t]
    \centering
    \includegraphics[width=0.496\linewidth]{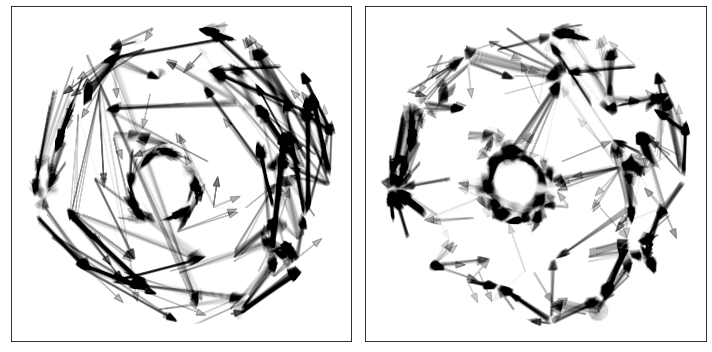}
    \includegraphics[width = 0.496\linewidth]{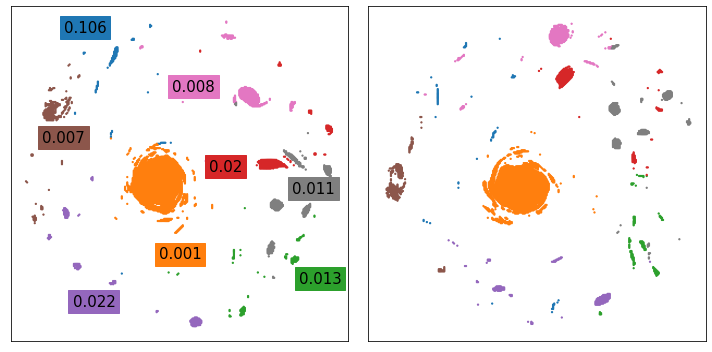}
    \caption{Left to right: GW transfer kernel, OT transfer kernel. The last two images show distinctly coloured partitions according to the spectral clustering of $K^{\GW}_{\eps,\kappa}$. The plotted numbers are the evaluations of the $\GW$ functional at $\pi_{\GW}$ restricted to the respective partitions.}
    \label{fig:vortex_ot_vs_fgw}
\end{figure}

We proceed to estimate the dynamics of the extracted patch. Similarly to our motivating example, we compare the OT and GW transfer operators. From $\omega_0$ and $\omega_1$ we extract the mm-spaces
$\XX = (X, d_X,\mu)$ and $\YY = (Y, d_Y,\nu)$, respectively. More precisely, $X$ and $Y$ are the sets of patch points in $\R^2$ where $\lvert \omega_0 \rvert \geq 600$ and $\lvert \omega_1 \rvert \geq 600$. Furthermore, $d_X$ and $d_Y$ are the normalized Euclidean metrics on $X$ and $Y$, respectively. Finally, $\mu$, $\nu$ are the (fully supported) probability measures proportional to the absolute value of the vorticity field. 
For our model we want to prohibit the transport between positive and negative vorticity.
To this end we label our data in the following way. Let $l_X,l_Y$ be the label function on $X,Y$ given by $0$, $1$ for negative, positive vorticity, respectively.
Additionally, due to possible dissipation of vorticity, we focus on unbalanced approaches for the estimation of the transfer operator.
We proceed to solve the entropic unbalanced OT problem between $\mu$ and $\nu$ with respect to the cost function  
$c(x,y) = d_{\text{E}}(x,y)^2 + d_\text{E}^2(l_X(x),l_Y(y))$,
regularization parameter $\eps = 0.0003$ and marginal relaxation parameter $\kappa = 0.1$.
This can be understood as the (entropic and unbalanced) Wasserstein distance with an additional penalty on transporting between distinctly signed vorticity. Thus we obtain a solution denoted $\hat{\pi}^{\W}_{\eps,\kappa}$. In the same way, let $\hat{\pi}^{\GW}_{\eps,\kappa}$ be a solution to the unbalanced, fused, entropic GW problem between $(\XX,l_X)$ and $(\YY,l_Y)$ and $\eps,\kappa$ as above.
On the left-hand side of Figure \ref{fig:vortex_ot_vs_fgw}, we illustrate the associated transfer operators $K^{\bullet}_{\eps,\kappa}$, $\bullet \in \{\W,\GW\}$. Similar to the previous examples, the OT transfer operator is not able to recover the underlying rotation. On the other hand, by favoring the preservation of intrinsic distances, the GW transport nicely reflects a counter-clockwise rotation.
Finally, we apply the discussed spectral clustering procedure, where we focus on $K^{\GW}_{\eps,\kappa}$. To obtain more than two coherent structures, we proceed in a nested manner. More precisely, applying the clustering procedure once yields two partitions of each mm-space $X_1,X_2$, $Y_1,Y_2$, respectively. Then we apply the procedure with respect to the associated (labelled) sub mm-spaces
\[
\biggl(\biggl(X_i, d_X \vert_{X_i}, \frac{\mu(\cdot \cap X_i)}{\mu(X_i)}\biggr),l_X\vert_{X_i} \biggr), 
\qquad \biggl(\biggl(Y_i, d_Y \vert_{Y_i}, \frac{\nu(\cdot \cap Y_i)}{\nu(Y_i)}\biggr),l_Y\vert_{Y_i} \biggr),
\]
and the restricted transfer operator $K^{\GW}_{\eps,\kappa} \vert_{X_i \times Y_i}$, $i=1,2$. This yields two sub-partitions per partition. 
We repeat this three times so that we obtain 8 partitions in total. The right-hand side of Figure \ref{fig:vortex_ot_vs_fgw} shows the mm-spaces $\XX,\YY$, where points of the same partition are coloured equally. Additionally, we evaluate the GW functional of $\hat{\pi}^{\GW}_{\eps,\kappa}$ restricted to the partitions 
as explained in Remark \ref{rem:1} and add the evaluation in the plot of $\XX$. 
As expected, the center vortex is clearly identified. 
Additionally, we are able to identify even smaller structures such as the coherent structures in brown, pink, grey and red. 
The orange partition attains the smallest $\GW$ evaluation by far. This is followed by brown and pink which represent smaller coherent structures and highlights their shape preservation under the transfer $K^{\GW}_{\eps,\kappa}$.


\section{Conclusions}
In this paper, we proposed a novel approach to estimate dynamical systems based on (unbalanced, fused) GW transport plans. Moreover, we demonstrated that the obtained transport plans can be leveraged for a spectral clustering procedure to extract coherent structures.
The resulting method is convenient as it can be quickly implemented by using out-of-the-box methods for GW and the singular value decomposition.
We verified its potency on three numerical examples.

As future work we leave a direct comparison with the method proposed in \cite{invOT}. The latter provides a numerically more appealing framework for obtaining transport plans which minimize the OT functional under additional invariances such as orthogonal transformations. Moreover, we are interested in applying our method on non-Euclidean data such as e.g.\ graphs.

\subsubsection{Acknowledgements} This work is supported by funds from the German Research Foundation (DFG) within the RTG
2433 DAEDALUS. The author thanks Jiahan Wang for fruitful discussions and his support regarding numerical implementations
as well as Gabriele Steidl for valuable discussions.
%
%
%
%

\end{document}